\documentclass{article}
\usepackage[utf8]{inputenc}

\usepackage[a4paper, total={6in, 8in}]{geometry}
\usepackage{amsmath,amssymb,amsthm}
\usepackage{hyperref}
\usepackage{graphicx}
\usepackage{color}
\usepackage{algorithm}
\usepackage{algpseudocode}
\usepackage{multirow}
\usepackage{authblk}
\usepackage{caption}
\usepackage{subcaption}

\usepackage{amsmath}
\usepackage{amssymb}
\usepackage{amsthm}
\usepackage{amssymb}
\usepackage{bbm}
\usepackage{hyperref}
\usepackage{color}

\newtheorem{theorem}{Theorem}
\newtheorem{lemma}{Lemma}
\newtheorem{remark}{Remark}
\newtheorem{cor}{Corollary}

\def \ER {Erd\H{o}s-R\'enyi }

\newcommand{\ceil}[1]{\lceil #1 \rceil}
\newcommand{\mb}[1]{\mathbb{#1}}

\newcommand{\brac}[1]{\left(#1\right)}
\newcommand{\cbrac}[1]{\left\{#1\right\}}
\newcommand{\sbrac}[1]{\left[#1\right]}
\newcommand{\indic}[1]{\mathbbm{1}{\brac{#1}}}
\newcommand{\expect}[2][]{\mathbb{E}_{#1}\sbrac{#2}}
\newcommand{\prob}[2][]{\mathbb{P}_{#1}\brac{#2}}

\newcommand{\abs}[1]{\lvert #1 \rvert}

\newcommand{\mf}[1]{\mathbf{#1}}

\title{Biased Consensus Dynamics on Regular Expander Graphs}
\author[$*$]{O. Deb}
\author[$\dag$]{A. Mukhopadhyay}
\affil[$*$]{Department of Physics, Jadavpur University}
\affil[$\dag$]{Department of Computer Science, University of Warwick}
\date{}

\begin{document}

\maketitle

\begin{abstract}
    Consensus protocols play an important role
in the study of distributed algorithms. In this paper, we study the effect of bias on two popular consensus protocols, namely, the {\em voter rule} and the {\em 2-choices rule}.
Classical versions of these rules assume that agents' behaviour do not depend on the state they are currently in. However, in many applications, such as social networks, there are often intrinsically superior states towards which agents exhibit some form of bias. In this paper, we are specifically interested in the case where the states or opinions of the agents are binary and come from the set $\{0,1\}$. We assume that agents with opinion $1$ update their opinion with a probability $q_1$ strictly less than the probability $q_0$ with which update occurs for agents with opinion $0$. We call opinion $1$ as the superior opinion and our interest is to study the conditions under which the network reaches consensus on this opinion. We assume that the agents are located on the vertices of a regular expander graph with $n$ vertices. 
This is a large class of graphs which includes both sparse (bounded degree) and dense (unbounded degree) graphs.
We show that for the voter rule, consensus is achieved on the superior opinion in $O(\log n)$ time with high probability if system starts with $\Omega(\log n)$ agents having the superior opinion.
This result is significantly different from the classical voter rule where consensus is achieved in $O(n)$ time and the probability of achieving consensus on any particular opinion is directly proportional to the initial number of agents with that opinion.  For the 2-choices rule, we show that consensus is achieved on the superior opinion in $O(\log n)$ time with high probability when the initial proportion of agents with the superior opinion is above a certain threshold. We explicitly characterise this threshold as a function of the strength of the bias and the spectral properties of the graph. We show that for the biased version of the 2-choice rule this threshold can be significantly less than that for the unbiased version of the same rule. Our techniques involve using sharp probabilistic bounds on the drift to characterise the Markovian dynamics of the system.
\end{abstract}

\section{Introduction}
\label{sec:intro}

Population protocols are studied extensively in distributed computing and statistical physics to probe how the group dynamics of a population of agents is affected by the local interaction rules which the agents follow to interact among themselves. 
A consensus protocol is  a set of interaction rules designed specifically to get all the agents agree on a common opinion.
It is of interest to study how fast and on what common opinion is consensus achieved.


The study of consensus protocols has received significant attention recently due to their applications in social networks. Here, the main interest is to find if social agents can learn better technologies or opinions through their interactions with neighbours. Another motivation 
comes from coordination problems in distributed systems where autonomous agents try to reach consensus by following simple rules to communicate with their neighbours. Such rules must be easy to implement, fault tolerant, and should lead to consensus quickly.

One of the simplest rules that have been studied in this context is the {\em voter rule}. Here, each agent in the network updates its opinion by copying the opinion of a randomly sampled neighbour. Due to its martingale nature and duality with coalescing random walks, the dynamics of voter rule has been studied extensively~\cite{holley1975ergodic,clifford1973model,cox1989coalescing,cooper2013coalescing,nakata1999probabilistic}. It has been shown that if all agents update synchronously in each time step, then, in a network with $n$ agents and two competing opinions, the voter rule leads to consensus in $O(n)$ steps for many classes of graphs including regular expanders. Furthermore, the probability of reaching consensus on a given opinion has been shown to be proportional to the volume (sum of degrees) of agents having that opinion initially.

Another rule which has been studied extensively in the literature is the 2-choice rule~\cite{cooper_two_choices,cooper_twochoice_gen,redner_two_choices}. Here, instead of sampling only one neighbour, an agent samples two neighbours uniformly at random and takes the majority opinion among the two sampled neighbours and the agent itself. This rule has been shown to lead to consensus much faster than the voter rule. Specifically, for regular expander graphs, this rule leads to consensus in $O(\log n)$ steps with high probability if the updates occur synchronously in each time step. However, in this case, the consensus is achieved with high probability on the opinion which is adopted by the majority of the population at the beginning. 

Most existing studies on the voter rule and the 2-choices rule with two competing opinions assume that the opinions are indistinguishable in the sense that an agent's behaviour does not depend on the opinion it currently possesses. However, in a social networking scenario, one opinion may represent a superior alternative such as a newer technology or a more effective drug for a certain disease. In such cases, agents already having the superior opinion are less likely to update their choices than agents with the alternative opinion. To model such scenarios, we introduce {\em bias} in the classical voter rule and in the 2-choices rule. Specifically, we represent the opinions as $0$ and $1$ and assume (without loss of generality) that agents with opinion $1$
update their opinion with a probability $q_1$ which is strictly less than the probability $q_0$ with which agents with opinion $0$ update their opinions. We assume that all updates occur synchronously in each round following either the voter rule or the 2-choices rule. Since agents with opinion $1$ are less likely to update in a given step than agents with opinion $0$, we call opinion $1$ the {\em superior opinion} and say that the population is biased towards this superior opinion.

This model of bias has been introduced recently in~\cite{arpan_JSP_20,mukhopadhyay2016binary}. Here, the authors study the dynamics under the biased voter rule and the biased 2-choices rule assuming that the underlying graph is complete. Their results indicate that consensus in the biased voter rule is achieved on the superior opinion with a probability which vanishes with increasing network size $n$ and the expected time to reach consensus grows as $O(\log n)$. This is significantly different from the unbiased voter rule where convergence to consensus is only linear in the network size and  the probability of achieving consensus on any given opinion remains bounded away from zero for all network sizes. For the biased 2-choices rule, their results indicate that consensus takes $O(\log n)$ time on average. Although this is similar to the unbiased version of 2-choice rule, the threshold required to achieve consensus on the superior opinion is now determined by the bias parameters $q_0$ and $q_1$. In particular, it is shown that if the initial proportion of agents with the superior opinion is above $q_1/(q_0+q_1)$ ($< 1/2$), then the network converges to consensus on the superior opinion with high probability. Thus, the network reaches a consensus on the superior opinion even if the superior opinion is initially held only by a minority of the population. The papers also studied the dynamics for more general classes of graphs numerically and conjectured that similar results hold for random regular graphs. However, an analytical study on general graphs was left as an open problem.

{\bf Our Contributions}: In this paper, we address this open problem. Specifically, we analyse the dynamics of the biased voter rule and the biased 2-choices rule assuming the underlying graph to belong to an expander family of regular graphs with increasing number of vertices. For simplicity, we shall call a member of such a family of graphs as a regular expander.
Regular expanders are a class of graphs whose degree for each vertex is the same (denoted by $d(n)$) and whose conductance, denoted by $\phi_n$, remains bounded away from zero for all network sizes, i.e., $\phi_n\geq \phi >0$ for all $n$ for some constant $\phi$. 
An equivalent way to define expanders uses the second largest (in absolute value) eigenvalue $\lambda_n$ of their scaled adjacency matrix. For an expander sequence of graphs, $\lambda_n \leq \lambda < 1$ for all $n$. Intuitively, the above properties imply that any large subset of nodes is connected to the rest of the network by a large number of edges; a larger value of $\phi$ or a larger spectral gap $1-\lambda$ implies a larger number of edges need to be removed to disconnect a subset of nodes with a given size.

Expanders are a large and important class of graphs often studied in the literature. It includes both sparse graphs
where $d(n)=O(1)$ and dense graphs where $d(n)=\omega(1)$.
Analysing consensus dynamics on such graphs is significantly harder than analysing the dynamics on complete graphs. This is because the Markov chain describing the evolution of the system has a higher dimensional state space and its transitions depend on the structural properties of the graph. The analysis is further complicated by the lack of symmetry in the biased model compared to the unbiased model. For example, the dynamics of the biased model no longer satisfies the martingale property and duality with coalescing random walks. To analyse such Markov chains we use the {\em drift} which is defined as the increase in the number of agents with the superior opinion in each step. By obtaining sharp probabilistic lower bounds on the drift, we prove the following results:

\begin{itemize}
    \item For the biased voter rule, we show that the network reaches consensus on the superior opinion in $O(\log n)$ steps with probability at least $1-O(\log n/n^\gamma)$, for some positive constant $\gamma$, if the initial number of agents with the superior opinion is $\Omega(\log n)$. Thus, even if the number of agents with the superior opinion is initially very small (logarithmic in the network size), with a high probability the network quickly reaches consensus on the superior opinion. Our analysis also shows that starting with a constant proportion of agents with opinion $1$ consensus occurs even faster with a higher probability. Thus, our results  significantly generalise the results obtained in~\cite{arpan_JSP_20,mukhopadhyay2016binary} for complete graphs.
    
    \item For the biased 2-choices rule, we show that if the expander family of graphs satisfies $\lambda_n^2 \leq q_0/(q_0+q_1)-c$ for each $n$ for some positive constant $c$ and if the initial proportion of agents with the superior opinion is larger than $q_1/(q_0+q_1)+\max\brac{\lambda_n^2,\sqrt{\log n/4n}}$, then consensus on the superior opinion is achieved in $O(\log n)$ time with probability at least $1-O(1/n)$. Note that for random regular graphs with degree $d(n)=\omega(1)$, we have $\lambda_n=o(1)$ with high probability. Hence, the condition $\lambda_n^2 \leq q_0/(q_0+q_1)-c$ is satisfied for all $n$ with high probability. In this case, we recover the results for complete graphs proved in~\cite{arpan_JSP_20,mukhopadhyay2016binary}. For random regular graphs with a constant degree $d$ the same condition  is satisfied with high probability for sufficiently large $d$ (for example, if $d > 4/(q_0/(q_0+q_1)-c)$). Therefore, our results are applicable to both dense and sparse graphs.
\end{itemize}

It is also worth mentioning here that our results hold even when an adversary is allowed to redistribute the opinions at each time step without changing the total number of agents having each opinion.

\section{Related Literature}

Opinion dynamical systems have been the subject of study for many years now. The simplest and the most well known model is the voter rule which was first studied independently in~\cite{holley1975ergodic} and~\cite{clifford1973model}. They discovered the martingale nature of the dynamics and its duality with coalescing random walks. Using these properties, the voter model has been studied for many important classes of graphs including random $d$-regular graphs~\cite{cooper2010multiple}, \ER  graphs~\cite{nakata1999probabilistic}, and regular lattices~\cite{cox1989coalescing}. Using the martingale property, it can be shown that for any connected graph the probability of reaching consensus on a specific opinion is proportional to the initial volume (sum of degrees) of nodes having that opinion~\cite{hassin2001distributed}. Furthermore, using the duality with coalescing random walks, it has been established in~\cite{cooper2010multiple} that consensus can be achieved on random $d$-regular graphs in $\Theta(n)$ time with high probability. 

The 2-choices rule studied in this paper was introduced as a faster alternative to the voter rule in~\cite{cooper_two_choices}. Here, the rule was studied for random regular graphs and expanders. It was shown that consensus can be achieved on the initial majority opinion in $O(\log n)$ time steps with high probability if the initial imbalance between the opinions is sufficiently high. These results were further generalised to in-homogeneous graphs without any regularity properties in~\cite{cooper_twochoice_gen}. Variants of the 2-choices rule have been also considered the physics literature~\cite{redner_two_choices,galam2002minority,chen2005consensus}. In these variants, groups of agents are formed in each time step, and all individuals in a group adopt the majority opinion within the group. These models have been shown to have similar convergence rates to consensus. A generalisation of 2-choices rule, where the updating agent samples $m$ agents from its neighbourhood and only changes its opinion if $d$ or more of the sampled agent differ from the updating agent, was analysed for complete graphs in~\cite{cruise2014probabilistic}. Similar rules have also been considered for random $d$-regular graphs in~\cite{abdullah2012consensus}. Their results imply that consensus can be achieved in $O(\log \log n)$ time when each agent samples at least five other agents in the network and takes the majority opinion.

The effect of bias on opinion dynamical models have been recent literature~\cite{mukhopadhyay2016binary,arpan_JSP_20,biased_info_sciences_2022,cruciani2021phase}. The model considered in this paper was first studied in~\cite{mukhopadhyay2016binary,arpan_JSP_20}. Here the authors study a continuous version of the model where each agent updates asynchronously at points of a Poisson process
associated with itself. Although this is different from the synchronous, discrete-time model considered in this paper, their results can be easily generalised to the setting considered in our paper. The analysis in these papers focused mainly on complete graphs whereas other graphs were studied numerically. The current paper generalises the results of these papers to regular expanders which include a large class of graphs including both sparse and dense graphs. A stronger form of bias has been considered in~\cite{biased_info_sciences_2022,cruciani2021phase}. Here, it is assumed that the agents can adopt the superior opinion with probability $\alpha$ irrespective of their current opinion and the opinions of their neighbours. Under such strong bias, consensus is only possible on the superior opinion. The convergence speed for this form bias has been analysed for the voter rule and the majority rule in~\cite{biased_info_sciences_2022}. The 2-choices rule under this form of bias has been analysed in~\cite{mukhopadhyay2022phase}. It is found that the convergence speed can vary drastically depending on the value of the parameter $\alpha$ and the initial proportion of agents with the superior opinion.

{\bf Organisation}: The rest of the paper is organised as follows. In Section~\ref{sec:model}, we introduce the models studied in the paper. In Section~\ref{sec:results}, we discuss the main results and their consequences. Here we also briefly discuss the outline of the proofs. The detailed analysis of the biased voter rule and the 2-choices rule are given in Sections~\ref{sec:analysis_voter} and~\ref{sec:analysis_twochoice}, respectively. Finally, we conclude the paper in Section~\ref{sec:conclusions}.

\section{Model}
\label{sec:model}

In this section, we introduce the model studied in this paper.
The model consists of $n$ interconnected agents; the interconnections are described by an undirected graph $G_n=(V_n,E_n)$ (with $\abs{V_n}=n$), where
each {\em vertex} or {\em node} $v \in V_n$ represents an  agent and each edge $(u,v) \in E_n$ represents a connection between two agents $u, v$.
 For each agent $u \in V_n$, we define $N_u = \cbrac{v:(u,v)\in E_n}$ to be the set of {\em neighbours} of $u$. 

Time is assumed to be discrete and at each discrete time step $t \in \mathbb{Z}_+$, each agent is assumed to
have an {\em opinion} in the set $\{0,1\}$.
We denote the opinion
of agent $u$ at time $t$ by $X_u^n(t) \in \cbrac{0,1}$.
The overall state of the network at any time $t \geq 0$ can be represented 
by the vector $\mf{X}^n(t)=(X_u^n(t), u \in V_n)$ of opinions of all the agents. 
The process $\mf X^n=(\mf X^n(t), t\geq 0)$, defined on the state space $\cbrac{0,1}^n$, describes the evolution 
of the system with time.

{\bf Bias}: We assume that at each discrete time step, each agent can potentially update its opinion and
the probability with which an agent does so depends on the agent's current opinion. 
More specifically,  if an agent currently has opinion $i \in \cbrac{0,1}$, then it updates its opinion with probability $q_i$ and chooses not to update its opinion with probability $1-q_i$.
To model bias, we assume that $q_0 > q_1$ , i.e., agents with opinion $1$ are less likely to update than agents with opinion $0$.
This is to model the fact that the agents which already have the superior alternative
are less likely to update their choice.

{\bf Update Rules}: If an agent updates its opinion at a given step,  then 
the update occurs according to a specific update rule.  We consider two popular update rules, namely, {\em the voter rule}
and the {\em 2-choices rule}.

{\em Voter rule}: Under the voter rule,
each agent
chooses another agent from its neighbourhood uniformly at random and copies its opinion in the next step.
Thus, for each agent $u \in V_n$ with opinion $i\in \{0,1\}$ at time $t$, the opinion at time $t+1$ is given by 
\begin{equation}
    \label{eq:voter_rule}
    X_{u}^n(t+1)=\begin{cases}
                        i & \text{ with probability }  1-q_i,\\
                        X_{N}^n(t) & \text{ with probability } q_i, 
                   \end{cases}
\end{equation}
where ${N} \in N_u$ denotes an  agent sampled uniformly at random from $N_u$.

{\em 2-choices Rule}: Under the 2-choices update rule,  each agent samples
two neighbours uniformly at random (with replacement) and updates to the majority opinion among the two sampled agents
and the agent itself. Thus for each agent $u\in V_n$ with opinion $i\in \{0,1\}$ at time $t$, the opinion at time $t+1$ is given by 
\begin{equation}
    \label{eq:2choices_rule}
    X_{u}^n(t+1)=\begin{cases}
                        i & \text{ with probability }  1-q_i,\\
                        M(u,N_1,N_2,t) & \text{ with probability } q_i, 
                   \end{cases}
\end{equation}
where $N_1$ and $N_2$ denotes two neighbours sampled uniformly at random from $N_u$, 
and $M(u, N_1,N_2,t)=\indic{X_{u}^n(t)+X_{N_1}^n(t)+X_{N_2}^n(t)\geq 2}$ 
denotes the majority opinion among two randomly sampled neighbours $N_1$ and $N_2$ of $u$ and the agent $u$ itself at time $t$. 




{\bf Graph Structure}: Throughout the paper, we assume the graph $G_n$ is $d_n$-regular for each $n$, i.e., $\abs{N_u}=d_n$ for all $u \in V_n$. Furthermore, we assume that the graph $G_n$ is a member of a sequence $\cbrac{G_n}_n$
which satisfies the following property:
\begin{equation}
    \phi_n \geq \phi >0, \forall n,
    \label{eq:conductance}
\end{equation}
where $\phi_n$ denotes the conductance of $G_n$ and $\phi >0$ is a constant independent of $n$.
We recall that the conductance $\phi_n$ of $G_n$ is defined as
\begin{align}
    \phi_n =\min_{S\subseteq V_n: \abs{S} \leq n/2}\frac{E(S,S^c)}{d_n\abs{S}}
    =\min_{S\subseteq V_n}\frac{E(S,S^c)}{d_n \min\brac{\abs{S},\abs{S^c}}},
    \label{conductance}
\end{align}
where for any subset $S \subseteq V_n$ of nodes of $G_n$, $E(S,S^c)=\abs{\cbrac{(u,v)\in E_n: u \in S, v\in S^c}}$ denotes the
number of edges connecting a node in $S$ to a node in $S^c=V-S$.
Intuitively, the property given by \eqref{eq:conductance} implies that, for every graph $G_n$ in the sequence, any large subset of nodes
is connected to the rest of the network by a large number of edges.
A higher value of the constant $\phi$ indicates better connectivity. 
A graph sequence $\cbrac{G_n}_n$ satisfying the above two properties 
is called a {\em $\phi$-expander family of regular graphs}. Expanders can also be characterised through their spectral properties.
Specifically, if $\lambda_n$ denotes the second largest (in absolute value) eigenvalue of the scaled adjacency matrix of $G_n$ and $\lambda_n \leq \lambda < 1$ for some fixed constant $\lambda \in(0,1)$ and for all $n$, then the graph sequence $(G_n)_n$ is called a {\em $\lambda$-expander family}.
The equivalence between the two definitions can be seen from Cheerger's inequality which states that $\sqrt{2(1-\lambda_n)} \geq \phi_n \geq (1-\lambda_n)/2$.

{\bf Absorbing states}: Under the assumptions stated above, the process $\mf X^n$ is a Markov chain on $\{0,1\}^n$ for each $n$.
Furthermore,
as a consequence of being an expander, each graph $G_n$ in the sequence $\{G_n\}_n$ is a connected graph (since $\phi_n \geq \phi > 0$ for each $n$).
This implies that the Markov chain $\mf X^n$
has only two absorbing states $\mf{0}$ and $\mf{1}$, corresponding to all agents having opinion $0$ and opinion $1$, respectively. 
Furthermore, since $G_n$ is connected, it is possible to reach any absorbing state from any non-absorbing state in a finite number of steps. 
Hence, with probability $1$, the chain $\mf X^n$ gets absorbed in either state $\mf 0$ or state $\mf 1$ in  finite time. 
We shall sometime refer to the absorption time as the {\em consensus time}. The main goal of our analysis   is to
find the conditions under which consensus is achieved on the preferred opinion $1$ and to characterise the consensus time as a function of the network size $n$, the network topology (described by $\phi$ or $\lambda$), the bias parameters ($q_0$ and $q_1$) and the initial state of the system given by $\mf X(0)$.

{\bf Notations}: Throughout our analysis we shall denote by $A_t=A(\mf X^n(t))$ (resp. $B_t=B(\mf X^n(t))$) the number of agents with opinion $1$ (resp. opinion $0$)
at time $t$. With a slight abuse of notation, we shall also use $A_t$ (resp. $B_t$) to denote the set of agents with opinion $1$ (resp. opinion $0$) at time $t$. We let $\Delta_{BA}(t)$ (resp. $\Delta_{AB}(t)$) denote the number of agents changing from opinion $0$ (resp. opinion $1$) to opinion $1$ (resp. opinion $0$) from time step $t$ to $t+1$. Furthermore, $\Delta(t)=A_{t+1}-A_t=B_t-B_{t+1}$ denotes the increase in the number of agents with opinion $1$ from time step $t$ to $t+1$. 
Unless mentioned otherwise, all expectations and probabilities are conditional on the current state of the process $\mf X^n$.

\section{Main Results}
\label{sec:results}

In this section we summarise our main findings and their consequences. We also discuss the main idea used to prove these results.
The technical details are given in subsequent sections.

For the biased voter model, our main result is the following theorem:

\begin{theorem}
\label{thm:voter}
Let $\cbrac{G_n}_n$ be a $\phi$-expander family of regular graphs.
Then, there exists a constant $\gamma \equiv \gamma(q_0,q_1,\phi) \in (0,1)$ such that under the biased voter rule 
if $A(0)=\Omega(\log n)$, then with probability
at least $1-O\brac{\log n / n^\gamma}$ consensus is achieved on superior opinion $1$ in $$O\brac{\frac{\log n}{\log \brac{{1+\frac{q_0-q_1}{2}\phi}}}}+O\brac{\frac{\log n}{\log\brac{\frac{1}{1-(q_0-q_1)\phi}}}}$$ time steps.
\end{theorem}

Note that in the above theorem each term appearing in the expression of consensus time grows as $O(\log n)$. However, they have different leading constants which depend on the expansion parameter $\phi$. They also correspond to two different phases of the dynamics (as explained in the analysis). This is why we keep them as separate terms in the theorem.
The theorem implies that in the biased voter model consensus is achieved with high probability in $O(\log n)$ steps on the superior opinion even when the initial number of agents having the superior opinion is only logarithmic in the network size. This is in sharp contrast to the classical voter rule where consensus takes $O(n)$ steps with high probability. Furthermore, the result is true for all positive values of the parameter $\phi$. This implies, that consensus occurs in $O(\log n)$ time even for expanders with poor expansion properties (i.e., with $\phi$ close to $0$). Also note that the speed of consensus depends on the bias parameters $q_0$ and $q_1$ through their difference. As long as this difference remains strictly positive, consensus occurs in logarithmic time.


For the 2-choices model, our main result is summarised in the following theorem:

\begin{theorem}
    \label{thm:2choice}
Let $\cbrac{G_n}_n$ be a $\lambda$-expander family of regular graphs with $\lambda^2 \leq q_0/(q_0+q_1)-c$
for some positive constant $c \in (0,q_0/(q_0+q_1))$.
Further, let $\gamma \leq c$ be a positive constant.  
If $(1-\gamma)\geq A(0)/n\geq q_1/(q_0+q_1)+\max\brac{\lambda_n^2,\sqrt{ \log n/4n}}$, then there exists a positive constant $\alpha=\alpha(q_0,q_1,\gamma,c) \in (0,1)$ such that with probability at least $1-O(1/n^\alpha)$ consensus is achieved on the superior opinion $1$ in $$O\left(\frac{\log n}{\log{\frac{1}{(1-q_1c(q_0/(q_0+q_1)-\gamma))}}}\right)$$ time steps. 
\end{theorem}

We first note that the condition $\lambda^2 \leq q_0/(q_0+q_1)-c$ implies that expanders should have sufficiently large spectral gap in order for consensus to occur in logarithmic time under the biased 2-choices rule. However, this condition is trivially satisfied by dense graphs where $d(n)=\omega(1)$ because for such graphs we have $\lambda_n=o(1)$. The condition of Theorem~\ref{thm:2choice} also holds with high probability for random $d$-regular graphs when $d$ is sufficiently high (e.g., $d\geq 4/(q_0/(q_0+q_1)-c)$. This follows from the fact that for such graphs $\lambda \leq (2\sqrt{d-1}+\epsilon)/d$ with high probability for any $\epsilon >0$. Therefore, the result of Theorem~\ref{thm:2choice} holds for a large class of graphs including both dense and sparse graphs. Furthermore, the theorem implies that for consensus to occur in logarithmic time the initial fraction of agents with the superior opinion must be above a certain threshold. Although this threshold depends on the spectral gap of the graph, for dense graphs (i.e., where $d_n=\omega(1)$) this threshold can be made arbitrarily close to $q_1/(q_0+q_1)$ by choosing $n$ sufficiently large since $\lambda_n=o(1)$ for such graphs. Note that $q_1/(q_0+q_1)$ was shown to be the asymptotic threshold value for complete graphs in~\cite{mukhopadhyay2016binary,arpan_JSP_20}. Hence, our result recovers the known results for complete graphs. 
According to Theorem~\ref{thm:2choice}, the consensus time decreases with increasing value of the constant $c$. Note that a higher value of $c$ is possible only when $\lambda$ is smaller, i.e., when the graph has good expansion properties.
This supports the natural intuition that the dynamics is expected to converge faster on graphs having better expansion properties.

The key to proving the above theorems is analysing the drift $\Delta(t)$ which captures the increase in the number of agents with opinion $1$ at each time step.
To have fast consensus on the superior opinion the drift should be non-negative in most time steps. To show this, we essentially follow three steps. The first step is to obtain a non-negative lower bound on the expected drift. In the second step, using concentration inequalities, we show that the drift in each step concentrates around its mean value with high probability. Thus, combining the first two steps, we  obtain a lower bound on the drift which holds  with high probability. The final step is to use this lower bound recursively in each step of the dynamics until consensus is reached. This is done by dividing the dynamics into two phases. The first phase brings the fraction of agents with superior opinion above a certain threshold value. The second phase further increases the number of agents with the superior opinion, ultimately resulting in a consensus.

\section{Analysis of the biased voter rule}
\label{sec:analysis_voter}

In this section, we prove the main result for the biased voter model. As mentioned before, the key is to analyse the drift $\Delta(t)$ of the Markov chain $\mf X^n$ at each step.
Throughout the section we assume that the graph sequence $\cbrac{G_n}_n$ is a $\phi$-expander family.
We first express the drift as a function of the current state of the network.

\begin{lemma}
\label{lem:voter_bound}
    For the biased voter model, the following hold.
    \begin{align}
        &\expect{\Delta_{BA}(t)} = q_0\frac{E(A_t,B_t)}{d_n}, \label{eq:EBA}\\
        &\expect{\Delta_{AB}(t)} = q_1\frac{E(A_t,B_t)}{d_n},\label{eq:EAB}\\
        &\expect{\Delta(t)} = (q_0-q_1)\frac{E(A_t,B_t)}{d_n}.
        \label{eq:voter_bound}
    \end{align}
\end{lemma}

\begin{proof}
For each agent $i \in B(t)$, let $Z_i(t)$ denote the Bernoulli random variable which takes the value $1$ if the agent $i$ changes to opinion $1$ at time $t+1$ and takes the value $0$, otherwise. Then, $\expect{\Delta_{BA}(t)}=\sum_{i \in B_t} \expect{Z_i(t)}$. Furthermore, under the biased voter rule,  $\expect{Z_i(t)}=\prob{Z_i(t)=1}=q_0 (d_i^A(t)/d_n)$, where $d_i^A(t)$ denotes the number of neighbours of agent $i$ in set $A$. Hence, $$\expect{\Delta_{BA}(t)}=q_0 \sum_{i \in B(t)} (d_i^A(t)/d_n)=q_0 E(A_t,B_t)/d_n,$$ 
where in the last equality we use the fact that $E(A_t,B_t)=\sum_{i \in B_t} d_i^A(t)$. This establishes~\eqref{eq:EBA}.
The proof of~\eqref{eq:EAB} is similar. 
Finally, combining~\eqref{eq:EBA}, \eqref{eq:EAB} 
and using the fact that $\Delta(t)=\Delta_{BA}(t)-\Delta_{AB}(t)$ we obtain~\eqref{eq:voter_bound}.
\end{proof}

The previous lemma shows that the expected drift is non-negative in each step. We now show that the drift is non-negative not just in expectation but also along the sample path of the process $\mf X^n$ with high probability.

\begin{lemma}
\label{lem:voter_prob}
For the biased voter model, the following holds for each $t\geq0$.
\begin{align}
    &\prob {\Delta(t)\leq \frac{q_0-q_1}{2} ~\phi ~ \min(A_t,B_t)}
   \leq 2e^{-\gamma\min(A_t,B_t)}
    \label{eq:voter_prob}
\end{align}
where $\gamma\equiv\gamma(q_0,q_1,\phi) \in (0,1)$ is a constant that depends only on $q_0$, $q_1$ and $\phi$. 
\end{lemma}

\begin{proof}
 We have 
\begin{align}
&\prob{\Delta(t)\leq \frac{q_0-q_1}{2} \phi \min(A_t,B_t)} 
\leq \prob{\Delta_{BA}(t)-\Delta_{AB}(t)\leq \frac{q_0-q_1}{2} ~ \frac{E(A_t,B_t)}{d_n}} \nonumber\\ 
&\leq \prob{\Delta_{BA}(t)\leq (1-\epsilon_1)q_0 \frac{E(A_t,B_t)}{d_n}}+\prob{\Delta_{AB}(t)\geq (1+\epsilon_2)q_1 \frac{E(A_t,B_t)}{d_n}}
\label{p1}
\end{align}
%
where $\epsilon_1=(q_0-q_1)/4q_0 \in (0,1)$ and $\epsilon_2=(q_0-q_1)/4q_1 >0$. In the first inequality, we have used the expander property of the graph and, in the second inequality, we have used union bound.

Now, we shall bound each term appearing on the RHS of~\eqref{p1}.
To do so we shall use the following Chernoff–Hoeffding inequalities: Let $Y=Y_1+Y_2+...+Y_N$ be the sum of the independent random variables $0\leq Y_i \leq 1$, $i= 1,2, . . . , N$. Then, for any $0\leq \epsilon_1 
\leq 1$ and any $\epsilon_2>0$ the following inequalities hold:
\begin{align}
\prob{Y\leq (1-\epsilon_1)\expect{Y}}  &\leq \exp\brac{-\epsilon_1^2 \expect{Y}/3},\label{ch1}\\
\prob{Y \geq (1+\epsilon_2) \expect{Y}} &\leq \Big(\frac{e^{\epsilon_2}}{(1+\epsilon_2)^{1+\epsilon_2}}\Big)^{\expect{Y}}\label{ch2}.
\end{align}
We have previously shown that both $\Delta_{BA}$ and $\Delta_{AB}$ can be expressed as a sum of independent Bernoulli random variables. Hence, we have

\begin{align}
&\prob{\Delta_{BA}(t)\leq (1-\epsilon_1)q_0\frac{E(A_t,B_t)}{d_n}} 
 \nonumber \\
& = \prob{\Delta_{BA}(t)\leq (1-q_0)\expect{\Delta_{BA}}} 
\leq \exp\brac{-{\frac{\epsilon_1^2}{3}~ \expect{\Delta_{BA}}}} \leq  \exp\brac{-{\frac{\epsilon_1^2}{3} q_0 \phi \min(A_t, B_t)}} \label{con1}
\end{align}
where $\epsilon_1=\frac{q_0-q_1}{4q_0} \in (0,1)$. In the above, the 
the first line follows from~\eqref{eq:EBA},
the second line follows from~\eqref{ch1}, and the last line  uses~\eqref{eq:EBA} and the expander property of the graph.

Similarly, using~\eqref{eq:EAB}, \eqref{ch2}, and the expander property of the graph we have 
\begin{align}
&\prob{\Delta_{AB}(t)\geq \Big(q_1 \frac{E(A_t,B_t)}{d_n}+ \frac{q_0-q_1}{4}\frac{E(A_t,B_t)}{d_n}\Big)} \nonumber \\
&=\prob{\Delta_{AB}(t)\geq\expect{\Delta_{AB}}\left(1+\frac{q_0-q_1}{4 q_1}\right)}\nonumber\\
& \leq  \Big(\frac{e^{\epsilon_2}}{(1+\epsilon_2)^{1+\epsilon_2}}\Big)^{\expect{\Delta_{AB}}}\nonumber\\
&=\exp\brac{-\epsilon_3 {\expect{\Delta_{AB}}}} \nonumber\\
& \leq 2 \exp\brac{-\epsilon_3 q_1\phi~ \min(A_t, B_t)},\label{con2}
\end{align}
where $\epsilon_2= \frac{q_0-q_1}{4 q_1} > 0$ and $\epsilon_3=-\log \frac{e^{\epsilon_2}}{(1+\epsilon_2)^{1+\epsilon_2}}>0$.
Using \eqref{con1}, \eqref{con2}, and \eqref{p1} we obatin the result of the lemma with $\gamma=\phi \min(\frac{\epsilon_1^2}{3} q_0, ~\epsilon_3 q_1)\in (0,1)$.
\end{proof}

\begin{lemma}
\label{lem:voter_prob1}
For any $i \leq n/2$ we have
\begin{equation}
\prob{B_{t+1}>n/2|B_t=i} \leq 2\exp\brac{- \frac{(q_0-q_1)^2\phi^2}{2}n}
\end{equation}
\end{lemma}

\begin{proof}
Let $C(t)$ be the set of nodes which have neighbours with the opposite opinion at time $t$. Furthermore, for $v \in B(t) \cap C(t)$,
let $Y_v(t)\in \{0,1\}$ denote the Bernoulli random variable which takes the value $1$ only if node $v$ changes its opinion from $0$ to $1$ at time $t$. Similarly, for $v \in A(t) \cap C(t)$, let $Y_v \in \{-1,0\}$ be a Bernoulli random variable which takes the value $-1$ only if node $v$ changes its opinion from $1$ to $0$ at time $t$. Then, we can express the drift as $\Delta(t)=\sum_{v \in C(t)} Y_v$. 

Now we use the following version of the Chernoff-Hoeffding inequality which states that if $Y_l$, $l=1,2,\ldots,N$ are independent random variables with $a_l \leq Y_l \leq b_l$ and $Y=\sum_{l=1}^N Y_l$ and $\eta>0$, then
\begin{equation}
    \prob{\abs{Y-\expect{Y}}\geq \eta}\leq 2\exp\brac{-2\eta^2/\sum_{l=1}^N (b_l-a_l)^2}.
\label{eq:ch_v2}
\end{equation}

Applying the above Chernoff-Hoeffding inequality with $Y=\Delta(t)$ yields
\begin{equation}
\prob{\abs{\Delta(t)-\expect{\Delta(t)}} \geq \eta}\leq 2\exp\brac{-2\eta^2/\abs{C(t)}}\leq 2\exp\brac{-2\eta^2/n},
\label{eq:ch_v3}
\end{equation}
where the last inequality follows because $\abs{C(t)}\leq n$.

Now we note that for any $i \leq n/2$ we have
\begin{align}
\prob{B_{t+1}>n/2|B_t=i} &= \prob{\Delta(t) \leq -(n/2-i)|B_t=i}\nonumber\\
&\leq \prob{\Delta(t)-\expect{\Delta(t)}\leq -(n/2-i)-(q_0-q_1)\phi i|B_t=i}\nonumber\\
&\leq \prob{\abs{\Delta(t)-\expect{\Delta(t)}}\geq n/2-(1-(q_0-q_1)\phi)i|B_t=i}\nonumber\\
&\leq 2\exp\brac{-2(n/2-(1-(q_0-q_1)\phi)i)^2/n}\nonumber\\
&\leq 2\exp\brac{- \frac{(q_0-q_1)^2\phi^2}{2}n},\nonumber
\end{align}
where the second line follows from~\eqref{eq:voter_bound}
and the fact that for $i \leq n/2$ we have $E(A_t,B_t)\geq d_n \phi i$; the fourth line follows from~\eqref{eq:ch_v3}; the last line follows 
since $(n/2-(1-(q_0-q_1)\phi)i)$ is minimised in the range $i \leq n/2$ when $i=n/2$.
\end{proof}

In the next part of our analysis, we divide
the evolution of the system, into two phases. Phase-I starts when opinion $1$ is the minority opinion and ends when the number of agents with opinion $1$ exceeds $n/2$ for the first time. Phase-II starts when opinion $1$ is already a majority and ends when a consensus is achieved on opinion $1$.
The next few lemmas characterise the time spent in each phase.
    
\begin{lemma}
\label{lem:voter_phase1}
Define $T_1=\ceil{\log \frac{n}{2A_0}/\log(1+\frac{q_0-q_1}{2}\phi)}$. Then, under the biased voter rule, starting with $A_0 < n/2$, opinion $1$ becomes the majority opinion for the first time in $T_1$ time steps with probability at least $1-2T_1e^{-\gamma A_0}$, 
where $\gamma$ is as defined in Lemma~\ref{lem:voter_prob}.
\end{lemma}
\begin{proof}
 Let $r=(1+\frac{q_0-q_1}{2} ~\phi)> 1$. We note that the definition of $T_1$ implies that $r^{T_1}A_0\geq \frac{n}{2}$.
Hence, if the events $F$ and $G$ are defined as $F=\cbrac{\exists t\in\{0,1,...,T_1-1\} \text{ s.t. } A_{t+1}< rA_t}$ and $G=\cbrac{A_{t} <\frac{n}{2} ~\forall t \in \{0,\ldots,T_1\}}$, then $F^c \subseteq G^c$ or $G\subseteq F$. Now, for each $\tau\in \cbrac{0,1,\ldots,T_1-1}$ define $F_\tau=\{A_{\tau+1}< rA_\tau \text{ and } A_{t+1} \geq rA_{t}, \forall t\in \{0,1,\ldots,\tau-1\}\}$ and $G_\tau=\cbrac{A_{t} <\frac{n}{2} ~\forall t \in \{0,1\ldots,\tau\}}$. Clearly, $F=\cup_{\tau=0}^{T_1-1}F_\tau$ and $G \subseteq G_\tau$. 
Hence, we have

\begin{align}
 \prob{G|A_0} =  \prob{G\cap F| A_0 }=\sum_{\tau=0}^{T_1-1}\prob{G \cap F_\tau | A_0 }\leq \sum_{\tau=0}^{T_1-1}\prob{G_\tau \cap F_\tau | A_0 }.\label{tmp13}
\end{align}
But we note that 
\begin{align}
\prob{G_\tau \cap F_\tau | A_0}
 &\leq \mb{P}\left(A_{\tau+1}\leq rA_{\tau}\bigg|A_0, 
rA_{t}< A_{t+1}< \frac{n}{2}~\forall t\in\{0,\ldots,\tau-1\}\right)\nonumber\\
&=\frac{\sum_{(\mf x(\tau),\ldots, \mf x(0))\in H}\prob{A_{\tau+1} \leq rA_\tau\vert \mf X^n(\tau)=\mf x(\tau)}\prob{\mf X^n(\tau)=\mf x(\tau),\ldots, \mf X^n(0)=\mf x(0)}}{\prob{(\mf X^n(\tau), \ldots,\mf X^n(0))\in H}}\label{tmp12}
\end{align}
where $H=\cbrac{A(\mf x(0))=A_0, rA(\mf x(t))< A(\mf x(t+1))< \frac{n}{2}~\forall t\in\{0,\ldots,\tau-1\}}$. Here, the equality in the second line follows from the Markov property of the process $\mf X^n$.  Now from Lemma~\ref{lem:voter_prob} we obtain
$$\prob{A_{\tau+1} \leq rA_\tau\vert \mf X^n(\tau)=\mf x(\tau)}\leq 2\exp\brac{-\gamma A(\mf x(\tau))}\leq 2\exp\brac{-\gamma A(\mf x(0))}=2\exp\brac{-\gamma A_0}$$ for each $\mf x(\tau)$ satisfying $(\mf x(\tau),\ldots,\mf x(0)) \in H$. 
Using this in~\eqref{tmp12} and using~\eqref{tmp13} we obtain 
$\prob{G|A_0}\leq 2T_1 \exp(-\gamma A_0)$ which completes the proof of the lemma.
\end{proof}

The above lemma implies that when $A_0$ is sufficiently large, then Phase-I ends in $T_1$ time with high probability.
In particular, for $A_0=\Omega(\log n)\leq n/2$ we obtain the following corollary. 
\begin{cor}
\label{cor:2}
If $A_0=\Omega(\log n) \leq n/2$, then Phase I completes in $O(\log n/\log(1+\frac{q_0-q_1}{2}\phi))$ time steps with probability at least $1-O({\log n}/{n^\gamma})$, where $\gamma$ is as defined in Lemma~\ref{lem:voter_prob}. 
\end{cor}

The next lemma characterises the time taken to complete Phase-II. 

\begin{lemma}
\label{lem:voter_phase2}
Define $T_2= \ceil{\frac{2\log n}{\log \frac{1}{1-(q_0-q_1)\phi}}}$. Then, under the biased voter rule, starting with $B_0 \leq n/2$, consensus is achieved on opinion $1$ in $T_2$ time steps with probability at least $1-O(1/n)$.
\end{lemma}

\begin{proof}
    If at an arbitrary time $t$, we have $B_t\leq \frac{n}{2}\leq A_t$, then, from \eqref{eq:voter_prob} we get, $$\prob{\Delta(t)\leq \frac{q_0-q_1}{2} ~\phi ~ B_t\bigg| B_t\leq n/2}= \prob{B_{t+1}\geq (1-\frac{q_0-q_1}{2} ~\phi) ~ B_t\bigg|B_t \leq n/2}\leq 2e^{-\gamma B_t}.$$
 
 %
 Now, we observe
 \begin{align}
 &\expect{B_{T_2}\vert B_0 \leq n/2}\nonumber\\ 
 &=\expect{B_{T_2} \cdot \left(\indic {B_{T_2-1}\leq \frac{n}{2}} + \indic {B_{T_2-1}> \frac{n}{2}}\right)\Bigg| B_0 \leq \frac{n}{2}}\nonumber\\
 &\leq \left(\sum_{i\leq \frac{n}{2}} \expect{B_{T_2}\cdot \indic {B_{T_2-1}= i}\vert B_0 \leq \frac{n}{2}}\right) 
 + n \expect{\indic {B_{T_2-1}> \frac{n}{2}}\vert B_0 \leq \frac{n}{2}}\nonumber\\
 &=\left(\sum_{i\leq \frac{n}{2}} \expect{B_{T_2}| B_{T_2-1}=i} \cdot \prob{B_{T_2-1}=i\vert B_0 \leq \frac{n}{2}}\right) 
 + n \prob{B_{T_2-1}> \frac{n}{2}\vert B_0 \leq \frac{n}{2}} \label{eq:temp1}
 \end{align}
We can further decompose the probability in the second term on the RHS of the last inequality as follows.

\begin{align}
   &\prob{B_{T_2-1}> \frac{n}{2}\vert B_0 \leq \frac{n}{2}} \nonumber\\
   =&~ \prob{B_{T_2-1}> \frac{n}{2},B_{T_2-2}\leq \frac{n}{2}\vert B_0 \leq \frac{n}{2} }+\prob{B_{T_2-1}> \frac{n}{2},B_{T_2-2}> \frac{n}{2} \vert B_0 \leq \frac{n}{2} }\nonumber\\
   \leq &~ \prob{B_{T_2-1}> \frac{n}{2},B_{T_2-2}\leq \frac{n}{2} \vert B_0 \leq\frac{n}{2}}+\prob{B_{T_2-2}> \frac{n}{2} \vert B_0 \leq \frac{n}{2}}\nonumber\\
   \leq &~ \sum_{t=0}^{T_2-2}\prob{B_{t+1}> \frac{n}{2},B_{t}\leq \frac{n}{2} \vert B_0 \leq \frac{n}{2}}\nonumber\\
   \leq &~ \sum_{t=0}^{T_2-2} \prob{B_{t+1}> \frac{n}{2}\vert B_0\leq \frac{n}{2}, B_{t}\leq \frac{n}{2}}.\label{eq:temp2}
\end{align}
Hence, using the last inequality in~\eqref{eq:temp1}, we obtain
\begin{align}
\expect{B_{T_2}\vert B_0 \leq n/2} &\leq \left(\sum_{i\leq \frac{n}{2}} \expect{B_{T_2}| B_{T_2-1}=i} \cdot \prob{B_{T_2-1}=i\vert B_0 \leq \frac{n}{2}}\right)\nonumber\\
&\hspace{1cm}+n\sum_{t=0}^{T_2-2} \prob{B_{t+1}> \frac{n}{2}\vert B_0\leq \frac{n}{2}, B_{t}\leq \frac{n}{2}}
\label{eq:temp5}
\end{align}
Since $\Delta(t)=B_t-B_{t+1}$, using~\eqref{eq:voter_bound} and the expander property we have
 \begin{equation}
 \expect{B_{t+1}| B_t\leq \frac{n}{2}} \leq r B_t, 
 \label{eq:temp3}
 \end{equation}
where $r=(1-(q_0-q_1))\phi < 1$. 

Furthermore, using the Markov property of $\mf X^n$ we have
\begin{align}
&\prob{B_{t+1}> \frac{n}{2}\vert B_0\leq \frac{n}{2}, B_{t}\leq \frac{n}{2}}\nonumber\\
&=\frac{\sum_{(\mf x(t),\mf x(0))\in H}\prob{B_{t+1}> \frac{n}{2}|\mf X^n(t)=\mf x(t)}\prob{\mf X^n(t)=\mf x(t), \mf X^n(0) =\mf x(0)}}{\prob{(\mf X^n(t),\mf X^n(0))\in H}}\nonumber\\
&\leq 2\exp\brac{-\frac{(q_0-q_1)^2\phi^2}{2}n} \label{eq:temp4}
\end{align}
where $H=\cbrac{(\mf x(t),\mf x(0)): B(\mf x(t)) \leq n/2, B(\mf x(0)) \leq n/2}$ and the last inequality follows from Lemma~\ref{lem:voter_prob1}.

Now, using \eqref{eq:temp3},\eqref{eq:temp4}, and \eqref{eq:temp5}, we get 
\begin{align}
 \expect{B_{T_2}\vert B_0 < \frac{n}{2}} &\leq   \sum_{i\leq \frac{n}{2}} r i ~\prob{B_{T_2-1}=i\vert B_0 < \frac{n}{2}}+ 2n(T_2-1)\exp\brac{- \frac{(q_0-q_1)^2\phi^2}{2}n}  \nonumber\\
 &\leq r \expect{B_{T_2-1}\vert B_0 < \frac{n}{2}}+ 2n(T_2-1)\exp\brac{- \frac{(q_0-q_1)^2\phi^2}{2}n}
 \label{exp1}
\end{align}
Using the above recursively we obtain
\begin{align}
   \expect{B_{T_2}\vert B_0 <\frac{n}{2}} 
   &\leq  r^{T_2} \expect{B_0}+2n  \exp\brac{- \frac{(q_0-q_1)^2\phi^2}{2}n}\sum_{k=1}^{T_2-1}(T_2-k)r^{k-1}\nonumber\\
   &\leq r^{T_2} ~\frac{n}{2}+2n  \exp\brac{- \frac{(q_0-q_1)^2\phi^2}{2}n} \frac{T_2}{1-r}
   \end{align}
Plugging in $T_2= \ceil{\frac{2 \log n}{\log \frac{1}{r}}}$ in the RHS we have
\begin{align}
  \expect{B_{T_2}\vert B_0 < \frac{n}{2}}\leq \frac{1}{2n}+   \frac{4n 
  \log n}{(1-r)\log \frac{1}{r}}\exp\brac{- \frac{(q_0-q_1)^2\phi^2}{2}n}= O\brac{\frac{1}{n}}
\end{align}
Hence, using the above and the Markov inequality we obtain,
\begin{align}
    \prob{B_{T_2}=0| B_0 < \frac{n}{2}}&= 1- \prob{B_{T_2}\geq 1| B_0 < \frac{n}{2}} \nonumber\\
    &\geq 1-\expect{B_{T_2} | B_0 < \frac{n}{2}} \nonumber\\
    &\geq 1- O\brac{\frac{1}{n}}, \nonumber
\end{align}
which proves the statement of the lemma.
\end{proof}

{\em Proof of Theorem~\ref{thm:voter}}: The statement of the theorem follows from Corollary~\ref{cor:2} and Lemma~\ref{lem:voter_phase2} using the Markov property of $\mf X^n$.\qed

\begin{remark}
From Lemma~\ref{lem:voter_phase1}, it follows that if $A_0=\ceil{np}$ for some constant $p >0$, then Phase-I completes in $O(1)$ time with probability at least $1-O(\exp(-n\gamma p))$. This implies that convergence to consensus on the superior opinion occurs even faster with a higher probability when a constant proportion of agents initially have the superior opinion in comparison to the case where only logarithmic number of agents have the superior opinion initially.
\end{remark}

\section{Analysis of the biased 2-choices rule}
\label{sec:analysis_twochoice}

In this section, we present the analysis of the biased 2-choices rule and the proof of Theorem~\ref{thm:2choice}. 
Throughout the analysis we shall assume that the sequence $\cbrac{G_n}_n$ is a $\lambda$-expander and shall use $\lambda_n$ to denote the second largest (in absolute value) eigenvalue of the scaled adjacency matrix of $G_n$.
As before, the first step in the analysis is to obtain a lower bound on the expected drift at each time step. This bound is more difficult to obtain in the 2-choices model than in the voter model since the drift has a non-linear dependence on the current state of the network. The following lemma characterises the expected drift in the 2-choices model.

\begin{lemma}
\label{lem:2choice_bound1}
    For the biased 2-choices model, we have
    \begin{equation}
        \expect{\Delta(t)}\geq B_t\left(q_1 (1-\lambda_n^2) \frac{A_t}{n} - \frac{q_1^2}{q_0+q_1} \right)
    \end{equation}
\end{lemma}

\begin{proof}
For each agent $i \in B(t)$, let $Z_i(t)$ denote the Bernoulli random variable which takes the value $1$ if the agent changes to opinion $1$ at time $t+1$ and takes the value $0$, otherwise. Then, $\expect{\Delta_{BA}(t)}=\sum_{i \in B_t} \expect{Z_i(t)}$. Furthermore, under the biased 2-choices rule,  $\expect{Z_i(t)}=\prob{Z_i(t)=1}=q_0 (d_i^A(t)/d_n)^2$, where $d_i^A(t)$ denotes the number of neighbours of agent $i$ in set $A$. Hence, $\expect{\Delta_{BA}(t)}=\sum_{i \in B(t)} q_0 (d_i^A(t)/d_n)^2$. Similarly, we have $\expect{\Delta_{AB}(t)}=\sum_{i \in A(t)} q_1 (d_i^B(t))/d_n)^2$. Hence, combining the above 
we obtain 
\begin{align}
\expect{\Delta(t)}
& = \sum_{i\in B_t} q_0\left(\frac{d^A_i(t)}{d_n}\right)^2-\sum_{i\in A_t} q_1\left(\frac{d^B_i(t)}{d_n}\right)^2,
\nonumber\\
& = q_1 \left(\sum_{i\in B_t} \left(\frac{d^A_i(t)}{d_n}\right)^2-\sum_{i\in A(t)}  \left(\frac{d^B_i(t)}{d_n}\right)^2\right)  +(q_0-q_1)\sum_{i\in B_t} \left(\frac{d^A_i(t)}{d_n}\right)^2 \label{eq:lower_split}
\end{align}
Thus, to obtain a lower bound on $\expect{\Delta(t)}$, it is sufficient to obtain a lower bound on each of the two terms on the RHS of the last equality. To bound the first term we use the following inequality from~\cite{cooper_twochoice_gen}:
\begin{align}
    &\sum_{i\in B_t} \left(\frac{d^A_i(t)}{d_n}\right)^2-\sum_{i\in A_t} \left(\frac{d^B_i(t)}{d_n}\right)^2 \geq B_t\left((1-\lambda_n^2)\frac{A_t}{n}-2 \theta_t (1-\theta_t)\right),
    \label{eqch4}
\end{align}
where $\theta_t=E(A_t,B_t)/d B_t$. To bound the second term we use Jensen's inequality to note that 
\begin{align} 
 \sum_{i\in B_t} \left(\frac{d^A_i}{d_n}\right)^2 &\geq \frac{1}{B_t d_n^2} (E(A_t, B_t))^2 = B_t\theta_t^2
\end{align}
Combining the above in~\eqref{eq:lower_split}, we obtain
\begin{align}
    \expect{\Delta(t)}
    & \geq q_1 B_t\left((1-\lambda_n^2)\frac{A_t}{n}-2 \theta_t (1-\theta_t)\right)+(q_0-q_1) B_t\theta_t^2\nonumber\\ 
    &=B_t\left(q_1 (1-\lambda_n^2) \frac{A_t}{n}-(2q_1 \theta_t - (q_0+q_1) \theta_t^2)\right).\nonumber
\end{align}
Since by definition $\theta_t$ always lies in the range $[0,1]$, the minimum value of the RHS of the last equality is obtained when $\theta_t=q_1/(q_0+q_1)$. Substituting this value in the last expression we obtain the lower bound as stated in the lemma.
\end{proof}

We note that in obtaining the lower bound in the above lemma we have not used 
the spectral properties of the graph $G_n$
We now use the spectral properties of the expanders to obtain a more refined lower bound on the expected drift which we shall use to analyse the dynamics under the 2-choices update rule. To state the result, we define $\epsilon(t)=(A_t-B_t)/n$ to be the imbalance between the two opinions at any time $t$ and $\epsilon'(t)=\epsilon(t)+(q_0-q_1)/(q_0+q_1)$ to be the shifted imbalance. In the following lemma, whose proof is given in the appendix, we obtain a more refined lower bound on the expected drift using the spectral properties of the graph.

\begin{lemma}
\label{lem:twochoice_bound2}
Assume that there exists some positive constant $c \in (0,q_0/(q_0+q_1))$ such that $\lambda_n^2 \leq q_0/(q_0+q_1)-c$ for all $n$. Then, for all $\epsilon' \geq 2\lambda_n^2$, we have
\begin{equation}
    \expect{\Delta(t)|\epsilon'(t)=\epsilon'}\geq \frac{B_t q_1}{2} c \epsilon'
    \label{eq:2choice_bound2}
\end{equation}
\end{lemma}

The above lemma implies that 
for any $\epsilon'\geq 2\lambda_n^2$, the following holds.
\begin{align}
\expect{B_{t+1}|\epsilon'(t)=\epsilon'}=\expect{B_t-\Delta(t)|\epsilon'(t)=\epsilon'}&\leq B_t\left(1-\frac{q_1 c}{2}\epsilon'\right).\nonumber
\end{align}
This implies that if the imbalance between the opinions is sufficiently high, then the number of agents with
opinion $0$ reduces in expectation at least by a factor. In the next lemma, we show that such reduction occurs not only in expectation but also along the sample path of the process $\mf X^n$ with high probability. 

\begin{lemma}
\label{lem:2choice_prob}
Assume that there exists a positive constant $c \in (0,q_0/(q_0+q_1))$ such that $\lambda_n^2 \leq q_0/(q_0+q_1)-c$ for all $n$.  Let $\gamma \leq c$ be a positive constant.
The following statements hold.

\begin{enumerate}
    \item For sufficiently large $n$ we have

\begin{equation}
    \prob{\Delta(t) \leq \frac{B_t q_1}{4} c \epsilon'\Bigg\vert  \epsilon'(t) =\epsilon'} \leq \frac{2}{n^\alpha},
    \label{eq:phase1_bound}
\end{equation}
where $\epsilon' \in \Gamma=\sbrac{\max\brac{2\lambda_n^2, \sqrt{\frac{\log n}{n}}},\frac{2q_0}{q_0+q_1}-2\gamma}$ and $\alpha=\frac{\gamma^2 q_1^2 c^2}{8} \in (0,1)$.

\item For any $i \leq n\gamma$, we have 

\begin{equation}
    \prob{B_{t+1} > n\gamma | B_t=i} \leq 2\exp(-n\beta),
\end{equation}
where $\beta=2\left(\frac{q_0}{q_0+q_1}-\gamma\right)^2 \gamma^2 q_1^2 c^2 > 0$.
\end{enumerate}
\end{lemma}

\begin{proof}
To prove the first statement of the lemma, note that $\gamma \leq c$ ensures that the interval $\Gamma$ is non-empty for sufficiently large $n$. 
Furthermore, $\epsilon'(t)\in \Gamma$ implies $B_t\geq n\gamma$, $\epsilon'(t)\geq 2\lambda_n^2$, and $n \epsilon'^2(t) \geq \log n$. 
As in the proof of Lemma~\ref{lem:voter_prob1}, we represent the drift $\Delta(t)$ as $\Delta(t)=\sum_{v \in C(t)} Y_v(t)$, where the set $C(t)$ and the random variables $Y_v, v\in C(t)$ have the same definitions as in the proof of Lemma~\ref{lem:voter_prob1}. We have

\begin{align}
    &\prob{\Delta(t)\leq \frac{B_t q_1}{4} c \epsilon'\Bigg\vert \epsilon'(t)=\epsilon'} \overset{(a)}{\leq} \prob{\Delta(t)\leq \frac{\expect{\Delta(t)|\epsilon'(t)=\epsilon'}}{2}\Bigg\vert \epsilon'(t)=\epsilon'}  \nonumber\\
    &\leq \prob{|\Delta(t)-\expect{\Delta(t)|\epsilon'(t)=\epsilon'}| \geq \frac{\expect{\Delta(t)|\epsilon'(t)=\epsilon'}}{2}\Bigg\vert \epsilon'(t)=\epsilon'} \nonumber \\
    &\overset{(b)}{\leq} 2 \exp\brac{-\frac{\mb{E}^2[\Delta(t)|\epsilon'(t)=\epsilon']}{2\abs{C(t)}}}
    \overset{(c)}{\leq} 2 \exp\brac{-\frac{B_t^2 q_1^2 c^2 \epsilon'^2(t)}{8n}}\overset{(d)}{\leq} \frac{2}{n^\alpha}\label{eq:temp_bound}
\end{align}
where (a) follows from~\eqref{eq:2choice_bound2}, (b) follows by applying Chernoff-Hoeffding inequality~\eqref{eq:ch_v3}, (c) follows from~\eqref{eq:2choice_bound2} and the fact that $\abs{C(t)} \leq n$, and (d) follows from the facts $B_t\geq n\gamma$ and $n\epsilon'^2(t)\geq \log n$.

To prove the second statement of the lemma, we note that
$B_t=i\leq n\gamma$ implies $\epsilon'(t)\geq 2q_0/(q_0+q_1)-2\gamma\geq 2\lambda_n^2$. Hence,
\begin{align}
\prob{B_{t+1}> n\gamma|B_t=i}&=\prob{\Delta(t)< -(n\gamma-i)|B_t=i}\nonumber\\
&=\prob{\Delta(t)-\expect{\Delta(t)|B_t=i}< -(n\gamma-i)-iq_1c\epsilon'/2|B_t=i}\nonumber\\
&\leq\prob{|\Delta(t)-\expect{\Delta(t)|B_t=i}|> n\gamma-(1-q_1c\epsilon'/2)i|B_t=i}\nonumber\\
&\leq 2\exp\brac{-2\frac{(n\gamma-(1-q_1c\epsilon'/2)i)^2}{n}}\nonumber\\
&\leq \exp(-n\beta)\nonumber
\end{align}
where the second line follows from~~\eqref{eq:2choice_bound2},
fourth line follows from~\eqref{eq:ch_v3}, and the last line follows from the facts that $\epsilon'\geq 2q_0/(q_0+q_1)-2\gamma$ and $i \leq n\gamma$.
\end{proof}

We divide the evolution of the system into Phase-I and Phase-II, where Phase I
brings the number of agents with opinion $0$  below $n\gamma$ and Phase-II further reduces it to $0$. The time taken in each phase is characterised in the following two lemmas. The proofs of these lemmas are similar to those of Lemma~\ref{lem:voter_phase1} and Lemma~\ref{lem:voter_phase2}. For completeness, we provide them in the appendix. 


\begin{lemma}
\label{lem:twochoice_phase1}
Assume that there exists a positive constant $c \in (0,q_0/(q_0+q_1))$ such that $\lambda_n^2 \leq q_0/(q_0+q_1)-c$ for all $n$ and
let $\gamma \leq c$ be a positive constant. 
If the system starts at an initial configuration satisfying $\epsilon'(0)\in \Gamma$, or, equivalently $$\gamma \leq \frac{B(0)}{n} \leq \frac{q_0}{q_0+q_1}-\max\brac{\lambda_n^2, \sqrt{\frac{\log n}{4n}}},$$ then, for sufficiently large $n$, the number of agents with opinion $0$ falls below $n\gamma$ for the first time in at most $T_1$ steps with probability at least $1-\frac{2T_1}{n^\alpha}$ where $T_1=\ceil{\log \brac{\frac{B_0}{n \gamma}}/\log(1-q_1 c \epsilon'(0)/4)^{-1}}$ and $\alpha$ is as defined in Lemma~\ref{lem:2choice_prob}.
\end{lemma}

\begin{remark}
We note that the initial condition $\epsilon'(0)\geq \max\brac{2\lambda_n^2, \sqrt{\log n/n}}$ in the lemma above can be satisfied for graphs with $\lambda_n=o(1)$ if $\epsilon'(0)$ is chosen to be an arbitrarily small positive constant and $n$ is sufficiently large. Hence, for sufficiently large dense graphs, the initial condition is satisfied by choosing $\epsilon'(0)$ to be any positive constant, or, equivalently, by choosing $B(0)/n$ to be any constant strictly less than $q_0/(q_0+q_1)$. This agrees with the results of~\cite{mukhopadhyay2016binary,arpan_JSP_20} on complete graphs. In general, choosing $\epsilon'(0) \geq 2\lambda^2$
satisfies the initial condition of the lemma for sufficiently large $n$ if the graph sequence is a $\lambda$-expander.
\end{remark}

\begin{remark}
We note that under the initial condition stated in the lemma $B_0/n < q_0/(q_0+q_1)$. Hence, $T_1=O(1)$ with probability at least $1-O(1/n^\alpha)$.
\end{remark}



\begin{lemma}
\label{lem:twochoice_phase2}
Assume that there exists a positive constant $c \in (0,q_0/(q_0+q_1))$ such that $\lambda_n^2 \leq q_0/(q_0+q_1)-c$ and let $\gamma \leq c$ be a positive constant.
Then, starting from $B_0 \leq n\gamma$ agents with opinion $0$, the network reaches consensus in $O(\log n/\log(1/(1-q_1c(q_0/(q_0+q_1)-\gamma))))$ time with probability at least $1-O(1/n)$.
\end{lemma}

{\em Proof of Theorem~\ref{thm:2choice}} The proof of the theorem now follows from Lemmas~\ref{lem:twochoice_phase1} and~\ref{lem:twochoice_phase2} using the Markov property of $\mf X^n$. \qed

\section{Conclusion and future directions}
\label{sec:conclusions}

In this paper, we have studied the dynamics of the voter rule and the 2-choices rule under the influence of bias assuming the underlying graph to be a member of an expander family of regular graphs. In our model, bias is introduced through the difference in update probabilities of agents having different opinions. We show that, under the biased voter rule, consensus is achieved on the superior opinion in logarithmic time with high probability even when the initial number of agents with the superior opinion is very small (logarithmic in the network size). For the biased 2-choices rule, our results indicate that consensus can be achieved in $O(\log n)$ steps with high probability on the superior opinion, provided that initial fraction of agents with the superior opinion is above a certain threshold. We characterise this threshold as function of the bias parameters and the spectral properties of the graph. 

The paper leaves several questions open for further investigation. One natural generalisation of the 2-choices rule studied in the paper is the $2k$-choices rule where an agent samples $2k$ agents from its neighbourhood and adopts the opinion of the majority among the sampled neighbours and the agents itself. It remains unclear how the number of agents sampled affects the dynamics when the agents are biased. Obtaining a sharp lower bound on the drift of the Markov chain will be challenging in this case. Another interesting generalisation to consider is the case where the underlying graph is not regular but satisfies the expander property. We believe similar bound on the consensus time should hold in this case but the conditions for achieving consensus would depend on the volume (sum of degrees) nodes having each opinion rather than the number of nodes having each opinion. The generalisation to the case with more than two opinions also remains as a challenging open problem.

\bibliographystyle{unsrt}
\bibliography{opinion.bib}

\appendix

\section{Proof of Lemma~\ref{lem:twochoice_bound2}}

From Lemma~\ref{lem:2choice_bound1} and the fact that 
$A_t=\frac{n}{2}(1+\epsilon(t))$ we obtain

\begin{align}
\expect{\Delta(t)|\epsilon'(t)=\epsilon'}&\geq B_t q_1 \left((1-\lambda_n^2) \frac{1+\epsilon(t)}{2} - \frac{q_1}{q_0+q_1} \right)\nonumber\\
&=\frac{B_t q_1}{2}\left(\epsilon(t)+ \frac{q_0-q_1}{q_0+q_1} - \lambda_n^2(1+\epsilon(t)) \right) \nonumber \\
&=\frac{B_tq_1}{2}\left(\epsilon'-\lambda_n^2 \left(\epsilon'+ \frac{2q_1}{q_0+q_1}\right) \right)\nonumber\\
&\geq \frac{B_tq_1}{2} c \epsilon',
\end{align}
where in the last line we use
the facts  $\epsilon'\geq 2\lambda_n^2$ and $\lambda_n^2\leq q_0/(q_0+q_1)-c$ which imply
$(1-\lambda_n^2)\epsilon'-\lambda_n^2 \frac{2q_1}{q_0+q_1}\geq c\epsilon'+\frac{q_1}{q_0+q_1}(\epsilon'-2\lambda_n^2)\geq c\epsilon'$.\qed

\section{Proof of Lemma~\ref{lem:twochoice_phase1}}

The proof of this lemma is similar to the proof of Lemma~\ref{lem:voter_phase1}. Define the event $F=\{\Delta(t)\leq B_t q_1 c\epsilon'(t)/4, \text{ for some } t \in \{0,1,\ldots,T_1-1\}\}$. Under $F^c$, we have $B_{t+1} \leq B_t\brac{1-q_1 c \epsilon'(t)/4}$ and $\epsilon'(t+1) \geq \epsilon'(t)$ for all $t \in \{0,1,\ldots,T_1-1\}$ which implies that $B_{T_1}\leq B_0\brac{1-q_1 c \epsilon'(0)/4}^{T_1}\leq B_0/r^{T_1}\leq n\gamma$, where the second inequality follows from the definition of $r$ and the last inequality follows from the definition of $T_1$. Hence, if $F^c$ occurs, then $B_{T_1}\leq n\gamma$. This implies that for $G$ defined as $G=\{B_t\geq n\gamma, \forall t\in \{0,1,,\ldots,T_1\}\}$, we have $G \subseteq F$.
Furthermore, we can partition $F$ as $F=\cup_{\tau=0}^{T_1-1} F_\tau$, where $F_{\tau}=\{\Delta(\tau)\leq \frac{B_\tau q_1 \epsilon'(\tau)}{4} \text{ and } \Delta(s) > \frac{B_s q_1 \epsilon'(s)}{4}, \forall s\in \{0,1,\tau-1\}\}$.
We have
\begin{align}
\prob{G|B_0}&= \prob{G\cap F|B_0}= \sum_{\tau=0}^{T_1-1}\prob{G\cap F_\tau|B_0}\leq \sum_{\tau=0}^{T_1-1}\prob{G_\tau\cap F_\tau|B_0}\label{tmp1}
\end{align}
where $G_\tau=\{B_t \geq n\gamma, \forall t \in \{0,1,\ldots,\tau\}\}$.
We further note that for sufficiently large $n$ we have
\begin{align}
\prob{G_\tau\cap F_\tau | B_0} &\leq \prob{\Delta(\tau)\leq \frac{B_\tau q_1 c\epsilon'(\tau)}{4} \Bigg\vert \Delta(s) > \frac{B_s q_1 c\epsilon'(s)}{4}, \forall s\in \{0,1,\tau-1\}, G_\tau, B_0} \nonumber\\
&\leq \frac{2}{n^\alpha},\label{tmp2}
\end{align}
where the inequality in the last line follows from the Markov property of $\mf X^n$, the fact that under $G_\tau$ and $\{\Delta(s) > \frac{B_s q_1 \epsilon'(s)}{4}, \forall s\in \{0,1,\tau-1\}\}$ we have $\epsilon'(\tau) \in \Gamma$, and the first statement of Lemma~\ref{lem:2choice_prob}. Finally, the statement of the lemma follows by combining~\eqref{tmp1} and~\eqref{tmp2}. \qed

\section{Proof of Lemma~\ref{lem:twochoice_phase2}}

Using the same steps as outlined in the proof of Lemma~\ref{lem:voter_phase2} we first obtain
\begin{align}
\expect{B_{T_2}|B_0 \leq n\gamma}&\leq \sum_{i \leq n\gamma}\expect{B_{T_2}|B_{T_2-1}=i}\cdot \prob{B_{T_2-1}=i|B_0 \leq n\gamma}\nonumber\\
&+n \sum_{t=0}^{T_2-2}\prob{B_{t+1}> n\gamma| B_t, \leq n\gamma, B_0\leq n\gamma}.    
\end{align}
But from the Markov property of $\mf X^n$ and second statement of Lemma~\ref{lem:2choice_prob} we have
\begin{align}
\prob{B_{t+1}> n\gamma| B_t \leq n\gamma, B_0 \leq n\gamma} \leq \exp(-n\beta)
\end{align}
Also, for $B_t \leq n\gamma$, we have $\epsilon'(t)=(n-2B_t)/n+(q_0-q_1)/(q_0+q_1)\geq  2q_0/(q_0+q_1)-2\gamma \geq 2\lambda_n^2$ since $\gamma \leq c$. Hence, from Lemma~\ref{lem:twochoice_bound2} we have
$\expect{B_{t+1}|B_t=i }\leq (1-\frac{q_1 c \epsilon'(t)}{2})i\leq ri$, for any $i \leq n\gamma$, where $r=(1-q_1 c (q_0/(q_0+q_1)-\gamma)) < 1$. Hence, $\sum_{i \leq n\gamma}\expect{B_{T_2}|B_{T_2-1}=i}\cdot \prob{B_{T_2-1}=i|B_0 \leq n\gamma}\leq r \expect{B_{T_2-1}|B_0\leq n\gamma}$.
Therefore, we have $$\expect{B_{T_2}|B_0 \leq n\gamma}\leq r \expect{B_{T_2-1}|B_0\leq n\gamma} + n(T_2-1)\exp\brac{-n\beta}.$$
Using the above inequality recursively we obtain
$$\expect{B_{T_2}|B_0 \leq n\gamma}\leq r^{T_2} B_0 + n\exp\brac{-n\beta}\sum_{k=1}^{T_2-1}(T_2-k)r^{k-1}\leq r^{T_2} n\gamma+n\exp\brac{-n\beta}(T_2/(1-r)).$$ Substituting $T_2=\ceil{2\log n/\log(1/r)}$ in the above inequality we obtain
$$\expect{B_{T_2}|B_0 \leq n\gamma}\leq \gamma/n+2n\log n\exp\brac{-n\beta}/((1-r)\log (1/r))=O(1/n).$$ Finally, using the above
and the Markov inequality the statement of the lemma follows. \qed

\end{document}